\documentclass[review]{elsarticle}

\usepackage[T1]{fontenc}
\usepackage{bm}
\usepackage{lineno}
\usepackage[colorlinks=true]{hyperref}
\usepackage{url}
\usepackage{amsmath,amssymb,amsthm}
\usepackage{mathpazo}
\usepackage{graphicx}
\usepackage{xcolor}
\usepackage{booktabs}
\usepackage{csvsimple}
\usepackage{siunitx}
\DeclareSIUnit{\fahrenheit}{\SIUnitSymbolDegree F}
\usepackage{fontawesome}
\usepackage{multicol,multirow}
\usepackage{csvsimple}

\usepackage{tikz}

\usepackage{subcaption}

\usepackage[margin=1in]{geometry}

\modulolinenumbers[1]

\journal{Elsevier}

\begin{document}
\begin{frontmatter}

\title{Demand Response Analogues for Residential\\ Loads in Natural Gas Networks}

\author{Fuyu Hu\fnref{fn1}}
\fntext[fn1]{School of International Affairs and Public Administration, Ocean University of China.}
\ead{hufuyu2010@gmail.com}

\author{Kaarthik Sundar\corref{cor1}\fnref{fn2}}
\cortext[cor1]{Corresponding author}
\fntext[fn2]{Information Systems and Modeling Group, Los Alamos National Laboratory.}
\ead{kaarthik@lanl.gov}

\author{Shriram Srinivasan\fnref{fn3}}
\fntext[fn3]{Applied Mathematics and Plasma Physics Group, Los Alamos National Laboratory.}
\ead{shrirams@lanl.gov}

\author{Russell Bent\fnref{fn3}}
\ead{rbent@lanl.gov}

\begin{abstract}
Demand response for electrical power networks is a mature field that has yielded numerous efficiency and resilience benefits like managing the peaks and valleys of electricity usage, reduction in peak electricity usage, to name a few. However, only recently has the study of the counterpart to demand response for electric power in natural gas started to receive similar levels of attention. Natural gas systems are increasingly operating at or near capacity, which challenges these systems to meet all the needs for gas, especially during severe winter weather. However, unlike demand response programs in electrical networks, demand response in gas networks cannot shift peak usage. Here, we develop analogues to demand response that can help improve the resilience of natural gas systems by reducing peak consumption and thereby limiting potential disruptions such events can cause. This paper develops a mathematical formulation to support a residential-level demand response analogue for natural gas based on current and anticipated smart thermostat technologies. The mathematical formulation takes the form of an optimal control problem (OCP) that leverages physics constraints to model temperature changes, balance equitable service, and optimize the gas consumption for collections of houses. On test problems, the formulation demonstrates significant benefits, including the ability to cut peak demand by 15\% while still ensuring equitable service to customers. 
\end{abstract}

\end{frontmatter}

\section{Introduction} \label{sec:intro}
One of the major challenges faced by modern energy systems is the problem of managing peak resource usage. During peak usage periods, energy systems operate at or near their capacity, which limits flexibility, increases green house gas emissions \cite{Kelley2018,EPRI2008}, reduces reliability \cite{Kevin2019}, decreases resilience, and, ultimately, drives costly infrastructure investments. To address this challenge, energy system operators are increasingly adopting and developing control solutions to help them better manage peak usage. One of the solutions adopted by the electric power industry is demand response (DR), which is used to shift energy usage away from peak consumption time periods. However, DR for other energy industries, such as natural gas, is relatively nascent.

DR in the electric power industry generally relies on two technologies: load control devices and market strategies. Market strategies are used to price energy and provide incentives to consumers to use energy at times desirable to the system operator \cite{Shariatzadeh2015}. Load control devices (energy storage, power-heating switch devices, etc.), are used to directly flatten energy demand by moving consumption from times with high utilization to times with low utilization. One example of load control in electric power DR programs is home cooling systems, where thermal inertia is used to pre-cool buildings prior to peak demand hours.

Within the electric power industry, 
there are numerous studies that have demonstrated the benefits of DR programs. For example, in 2009, the United States Federal Energy Regulatory Commission (FERC) announced a national DR action plan for maximizing the adoption of DR. The plan included a national communication program, the development of analytical tools, the outlining of regulations, and examples of model contracts \cite{FERC2009,FERC2010}. 
Subsequent reports estimated that adoption of DR has saved approximately 2,133 \si{\mega\watt} of peak demand in wholesale markets (2018) and approximately 30,000 \si{\mega\watt} of peak demand in retail markets (2017). 
While natural gas is an important energy resource \cite{ng-increase}, it has received little DR attention, as compared to electric power, and much of the existing work has focused on how electric DR can support natural gas by, for example, indirectly reducing natural gas used to produce electricity.
Unfortunately, winter months see a considerable increase in non-electric natural gas usage (e.g., home heating) for which electric power DR offers no relief.
Moreover, extreme weather events, like the 2021 cold weather snap in Texas, stress the natural gas infrastructure in ways that an electric power DR implementation cannot support.
Hence, there is an \emph{important opportunity for gas specific DR programs, e.g., a DR program focused on residential home heating}. 


A DR program focused on residential home heating differs from similar programs for electric power DR for several reasons. 
In an electric power DR program, the operator introduces variable pricing to encourage greater use during off-peak hours; the consumer weighs the options and determines if it is viable to shift usage to periods with lower pricing, and usually, energy usage is shifted away from peak time periods. However, when considering a natural gas DR program for residential home heating, it is apparent that the aims and ideas of an electric power DR do not directly apply. For instance, variable pricing of gas may exacerbate social inequities for a lifeline use of natural gas when it is most needed (e.g., cold temperatures). 

\emph{This observation reveals the question at the heart of this paper--if demand response, as usually understood, does not translate to gas DR, is there an analogue that can be introduced?
\footnote{In the subsequent discussion, the terms natural gas DR or gas DR are used with the tacit understanding that the terms do not correspond in a natural way to their counterpart in electric power DR.} While an operator cannot shift peak demand, it is still useful to reduce total and peak consumption, either at the level of an individual household, or for an ensemble of houses, especially in extreme winter weather conditions such as a polar vortex.  For consumers, it requires agreement to allow control of heating equipment to be guided by the forecast of the ambient temperature and current internal temperature, along with  a willingness to give up some amount of thermal comfort (in the form of deviation from set-point temperatures) for monetary incentives. A reduction in consumption could be the difference between total collapse of the distribution network (as happened in Texas recently \cite{texas2021}) and keeping that network functioning through challenging conditions with an equitable distribution of gas supply.}

There are many challenges associated with natural gas DR \cite{FOA2519} and this paper focuses on \emph{addressing the challenge associated with the formulation of mathematical models based on physics and thermodynamics for controlling natural gas consumption}. To the best of our knowledge, this is the first effort to formulate and demonstrate the workings of such a model to support natural gas DR programs. The model assumes current and anticipated capabilities of smart thermostats and is focused on smaller consumers, who are a significant source of DR capability in electric power systems, but are yet to be studied in natural gas systems. Since home heating systems are the major source of natural gas consumption in residential households, the model developed in this paper optimizes the gas consumption (home heating) of a collection of houses to reduce peak gas demand during normal winter operation, as well as during extreme temperature events, such as those induced by severe weather events.

While a DR program has various aspects (incentive programs, objectives, thermostat controller options, etc.),  all DR strategies are capable of being formulated as an optimal control problem (OCP) from the perspective of  mathematical modeling. The OCP we develop is a mixed-integer optimization problem  coupled with ordinary differential equations (ODEs) that model the dynamics of house temperature (heat gain from gas-fired furnaces and heat loss to the surroundings due to cooler external temperatures). The details of the OCPs are driven by two factors: (i) the implementation of the control policy on smart thermostats (centralized or decentralized) and (ii) the objective of the OCP (e.g., maintaining household temperatures as close to desired set points as possible while keeping peak consumption below a specified threshold, or reducing the overall gas consumption while keeping temperature deviations to a pre-specified limit). Overall, the goal of all the OCPs presented in this article is to assess (qualitatively and quantitatively) the merits  of a DR policy through the viewpoint of both the natural gas system operator and residential customers. The models also serve as an  important first step towards identifying requirements for DR technologies and market structures. We expect the methods developed in this paper to be used by stakeholders as a baseline for identifying potential benefits and determining the underlying requirements for a successful implementation of DR in a natural gas setting. In summary, while the crux of this article is the mathematical formulation of OCPs that can model gas DR for residential use, the key contributions of this paper include:

\begin{itemize}

\item Formulation of an Ordinary Differential Equation (ODE) model for a home heating system inspired by dynamic models for Thermostatically Controlled Loads (TCLs) and leverage these to formulate OCPs that optimize natural gas load scheduling at a household level.

\item Illustration of the use of a receding-horizon based algorithm to solve the resulting OCPs when the discretized problem becomes large.

\item Presentation of case studies which highlight the aggregation benefits and the degree of thermostat flexibility required to achieve a certain level of peak reduction. 

\end{itemize}

The remainder of this paper is organized as follows: Section \ref{sec:background} provides a literature review and relevant background for DR. Section \ref{sec:setup} discusses a high level description of the problem. In Section \ref{sec:overview}, we present a baseline model to reflect current practice and formulate the physics-based, decentralized and centralized natural gas DR models and  describe the solution method in Section \ref{sec:algo}. Finally, Section \ref{sec:case-study} presents the case studies and simulation results and Section \ref{sec:future} summarizes the conclusions and offers future perspectives.

\section{Background and Literature Review}\label{sec:background}

Optimization is the core computational technology used for investigating and controlling DR. Within the literature, both deterministic and stochastic optimization methods have been developed. Deterministic models assume that all information is known by the users and the operators \cite{Goudarzi2011, Barbato2014}. Stochastic models introduce uncertainty associated with supply, demand, price, customer behavior, environmental variables, etc. \cite{Massrur2018, Wang2016, Bai2016, Zhang2016, Zhang2020, He2019, Ju2019, Su2017}. The objective functions consider economic factors like cost savings \cite{Bradley2013, Qadrdan2017}, operation profits \cite{Ju2019, Sheikhi2015, Cui2016, Zhou2018}, engineering, reliability \cite{WangF2017}, risk minimization \cite{Ju2019}, greenhouse gas emission reductions \cite{Su2020, Su2017}, social welfare \cite{Sheikhi2015, He2019} and user comfort \cite{Wang2019}. Other lines of work include optimal control of ensembles of Thermostatically Controlled Loads (TCLs) \cite{Halder2019,Bhattacharya2017,Bhattacharya2018} which are very closely related to the models and DR formulations that are presented in this paper. Ensemble-based optimal control problems reduce electricity consumption from TCLs (typically air-conditioners) during the summer months by utilizing thermal inertia. In this line of work \cite{Halder2019}, a simple Newtonian thermal model (an RC circuit model) for heat exchange between the surroundings and a house is utilized. \emph{We use similar thermal models in this article to develop a DR program in a natural-gas setting to cater to the heating needs of households in winter months}. There are two fundamental differences between conventional (electric power based) TCL and the natural gas furnace based control considered here. First, TCLs operate in steady-state, while gas furnace operation is inherently transient. Second, TCLs have almost instantaneous response, while the gas furnaces operate on much longer time-scales, thus ruling out methods based on thermal-inertia that were successful for TCLs.

While DR response has been a part of the electric power industry's technology adoption conversations for at least a decade, it is only recently that such technologies have been discussed in the context of other energy industries, such as natural gas. This interest is largely driven by two recent emerging drivers.
First, historically low prices have greatly expanded natural gas's utilization and placed enormous strains on existing capacity \cite{babula2014cold}. At the same time, economic, political, and social factors have limited the ability of the industry to add capacity \cite{post2020}. As a result, DR is being proposed as a technology for improving the utilization of existing capacity \cite{House2019}. Second, extreme events like polar vortices and  the recent 2021 cold snap in Texas have created sudden spikes in the natural gas consumption for home and commercial heating that contributed to loss of load to millions of homes \cite{texas2021}. Such events can create supply and capacity shortages that interrupt gas deliveries. During recent events, states like Minnesota have employed an extreme form of DR that asked customers to decrease thermostat set-point temperatures to  \SI{60}{\fahrenheit} to counter such demand spikes \cite{polar-vortex}. \emph{Despite the potential for DR to support the natural gas industry in both these contexts, work and approaches for modeling and controlling DR in natural gas systems has been limited (unlike electric power) but is expected to increase dramatically in the next several years \cite{FOA2519}.}

There is also some research on natural gas as a supplemental energy resource for electric power DR.
\cite{Cui2016} developed a DR market model where natural gas based generation provides a supportive role for electric power DR.  \cite{Qadrdan2017} developed a model of the combined electric and natural gas system in Great Britain and used it to study the influence of DR in combined system modeling. \cite{Su2017} used electric power DR technologies, like pumped storage and incentives, to improve the operational performance of coupled natural gas and electric systems. More recently, \cite{Zhang2018} developed a pricing strategy to maximize the social benefits of DR in coupled gas-grid systems. 
Finally, \cite{Abahussain2013} described a price-based DR for optimally scheduling deferable gas loads (power plants) to maximize profits and \cite{Rudkevich2017} built on this approach to develop a combined pricing mechanism to schedule production and consumption of natural gas systems to maximize social welfare. 
In all these papers, the focus was largely on how DR in power systems could impact the operations of joint gas-electric systems with \emph{limited attention devoted to how to use gas assets in DR}.

\section{Demand Response Setup and Assumptions} \label{sec:setup}
\begin{figure*}[htb]
    \centering 
    \includegraphics[scale=0.2]{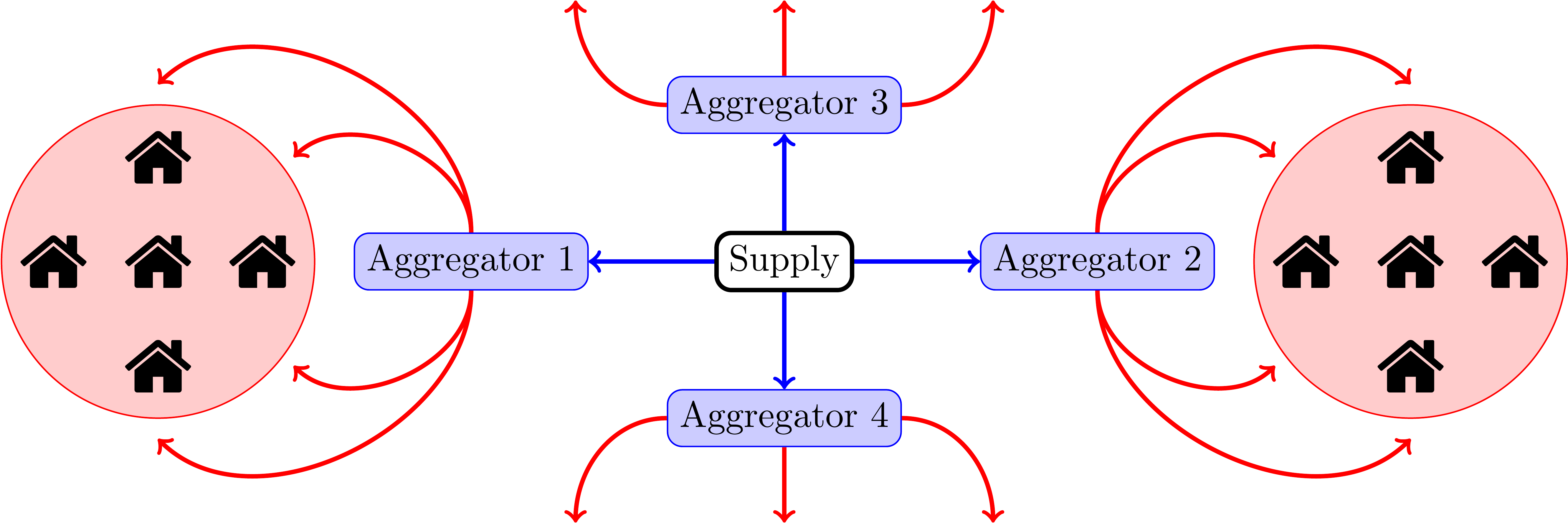}
    \caption{An illustrative diagram of a DR system. Aggregators implement a DR program for a set of houses that it manages. Here, an aggregator may refer to a utility company or other third party in a particular geographic location.}
    \label{fig:dr_schematic}
\end{figure*}

Figure \ref{fig:dr_schematic} shows a conceptual diagram of consumer level DR programs. 
While this is not the only way DR programs are structured, it is a common method and this structure is used to inform the OCP and solution methods developed in this paper.
In this figure, the operator of the DR program seeks to reduce energy consumption during peak periods. Consumers (households) are given incentives to participate in the DR program to help meet the goals of the system's operator. Load aggregators (blue boxes in Fig. \ref{fig:dr_schematic}) are sometimes used to recruit participants, serve as an interface between the individual consumers and markets, and are used to coordinate consumer responses to DR requests.  

A key component of the gas DR of this paper is ensuring relatively equal distribution of adverse conditions to participants. In this context, an adverse condition is thermal discomfort (modeled as average deviation from temperature set points or the maximum deviation from temperature set points) which has detrimental effects on human health and performance \cite{Lan2011,Jevons2016,LAN20111057}. In the OCP, the requirement of equitable distribution for customer satisfaction is enforced by including a social welfare objective function which indirectly minimizes the thermal discomfort. Another component of DR is participant incentives as participants need incentives (typically monetary) to give up temperature comfort. 
\emph{This paper uses the OCP to identify key parameters 
(maximum temperature deviation and number of participants) to inform the choice of incentive levels but does not directly address the challenge of implementing incentives.}  
The final component of DR is 
constraints that model the governing equations describing the dynamics and control of the system. 
Whereas electric power based DR is well developed, \emph{natural gas DR is relatively new and these governing equations have yet to be fully established}.
In this article, we utilize thermodynamic models of heat flow developed for electric power TCLs \cite{Halder2019,Bhattacharya2018} and contribute a natural gas DR formulation for these dynamic models.

To construct the OCP, we make the following assumptions about DR. These assumptions allow us to simplify the mathematical model, and investigate  the worst-case guarantee that can be provided for the benefits a natural gas DR program can support. The assumptions include:

\begin{enumerate}
    \item The natural gas utility operator aggregates the natural gas demand from all households in a locality. From here on, we use the term natural gas utility operator and aggregator interchangeably. This assumption is a mere statement of fact on the functioning of natural gas utility operators. 
    \item The only source of natural gas consumption for the households is home heating. This is justified because home heating systems are the major source of flexible natural gas consumption in most households during the winter months. This assumption enables us to ignore other minor sources of natural gas consumption at a household level like natural gas cook tops.
    \item Each participating household in the DR program has a smart thermostat. In this article, we consider two common thermostats (``Type 1'' and ``Type 2'') and their use in an OCP. The ``Type 1'' thermostat is capable of (a) obtaining a temperature forecast for the future, (b) allowing the customer to establish a temperature set-point for the thermostat, (c) allowing the customer to choose a parameter that adjusts the priority of the thermostat between maintaining a set-point temperature to saving on gas consumption, and (d) has minimal computing capability to solve optimization problems. The use of a ``Type 1'' thermostat permits a decentralized control policy as no information is shared between households or between the household and the aggregator. 
    
    The ``Type 2'' thermostat is capable of (a) sending the current temperature of the house and the set-point temperature of the thermostat to the aggregator and (b) receiving on-off control policy from the aggregator that it implements. These controllers support centralized control policies. Unlike the ``Type 1'' thermostat, the ``Type 2'' thermostat does not perform any computation and all of the computation is performed by the aggregator.  In both cases, \emph{the underlying OCP is largely the same and relies on the same algorithm developed in this paper.}
    \item Effects of the distribution pipeline network used to deliver gas to the households are saved for future work. 
\end{enumerate}

For this viability study, we argue that if one does not quantitatively or qualitatively find any benefit to using the OCP model for DR under these assumptions, then any DR formulation with less strict assumptions will only do worse and be unviable. Furthermore, the effectiveness of utilizing the proposed DR formulations are quantitatively corroborated using a combination of (i) the comfort for each participant examined in terms of statistics of temperature deviation from the set-point temperature and (ii) the overall reduction in the amount of gas consumed by the participant. On the other hand, benefits to the aggregator are quantified using the overall reduction in the gas consumption achieved when more participants sign-up for the DR programs.

\section{DR optimal control problems} \label{sec:overview}
Before presenting Optimal Control Problem (OCP) formulations for the natural gas DR problem, we introduce some notation. A brief outline of the notation is provided in the nomenclature section below while a detailed description of the entries is provided in context as and when they are used to formulate the OCPs. 

\subsection{Nomenclature}\label{subsec:nomenclature} 
\noindent \textit{Sets and parameters:} \\
$\mathcal H$ - set of houses indexed by $h$ \\
$\bm \Theta(t)$ - forecast of ambient temperature \\
$\bm \alpha_h, \bm \beta_h$ - thermal coefficients of house $h$ \\
$\bm Q_h$ - burn rate of the gas furnace in the house $h$ \\
$\bm \theta_h^0$ - initial temperature of house $h$ \\
$\bar{\bm \theta}_h$ - set-point temperature of house $h$ \\
$\bm \gamma$ - load curtailment factor \\ 
$T$ - optimization time horizon \\
$\bm \lambda_h$ - thermostat parameter that controls its mode of operation \\
$\bm D$ - max. total gas consumed by all the households in $\mathcal H$ \\ 
$\bm V_h$ - volume of house $h$ \\ 
$\bm \kappa_h$ - heat transfer coefficient of the walls of the house $h$ \\
$\bm A_h$ - surface area of the walls of the house exposed to the outside temperature \\ 
$\bm E^g$ - specific heating value of natural gas \\ 
$\bm C^a$ - isochoric specific heat of air \\ 
$\bm \rho^a$ - density of air \\
\noindent \textit{State variables:} \\
$\theta_h(t)$ - temperature of house $h$\\
\noindent \textit{Control variables:} \\
$q_h(t)$ - on-off variable for the natural gas supply of house $h$\\
\noindent \textit{Other variables:} \\
$\Delta_h$ - maximum deviation of the house temperature from the set-point temperature of the house \\
$G_h$ - total gas consumption for the house from $[0, T]$

We next present the thermodynamics of heating inside a single house. For ease of readability, we present the explanation of the notation in detail. 

\subsection{Thermodynamics of heating} \label{subsec:dynamics}
In a locality where a DR program is implemented, we denote by $\mathcal H$ the set of houses, and index  each house by $h$. We let $\bm \Theta(t)$ be the given time-varying forecast of the ambient/outside temperature. We let $\bm \theta_h^0$ denote the initial temperature inside house $h$. Each house, $h$, has a gas furnace that burns fuel at the rate $\bm Q_h$  (\si{\kilogram\per\second}). For each house $h \in \mathcal H$, we let $\{\bm \alpha_h, \bm \beta_h\}$ denote its thermal coefficients. Also associated with each house $h$ are the following state and control variables: $\theta_h(t)$ - current temperature of the house (\si{\kelvin}) and $q_h(t) \in \{0,  1\}$ - a binary control variable that controls supply  of natural gas. For each house, $h \in \mathcal H$, we denote the differential algebraic equations (DAEs) that govern the thermodynamics of heating by the following equations:
\begin{subequations}
\begin{flalign}
    & \dot{\theta}_h(t) = -\bm \alpha_h \cdot (\theta_h(t) - \bm \Theta(t)) + \bm \beta_h \cdot \bm Q_h \cdot q_h(t) \quad \forall t \label{eq:ode-h}, \\
    & \theta_h(t=0) = \bm \theta_h^0, \text{ and } \label{eq:ic-h}\\
    & q_h(t) \in \{0, 1\}. \quad \forall t \label{eq:control-h}
\end{flalign}
\label{eq:dynamics}
\end{subequations}

\noindent
Eqs. \eqref{eq:ode-h} is a  Newtonian equation that captures the dynamics of heat flow within the house $h$. A similar equation is used to model TCLs in the electric power DR literature \cite{Halder2019}. In Eq. \eqref{eq:ode-h}, $\dot{\theta}_h(t)$ is used to indicate the time derivative of the temperature of the house. The algorithm to compute the thermal coefficients $\{\bm \alpha_h, \bm \beta_h\}$ for each house $h \in \mathcal H$ is provided in the Appendix. Eqs. \eqref{eq:ic-h} -- \eqref{eq:control-h} define the initial conditions for $\theta_h(t)$ and the binary restrictions on the on-off control variable $q_h(t)$, respectively. 
While this model incorporates the effect of the ambient temperature, it does not take into account the effect of solar radiation or the behaviour of the occupants within the house. 
We use the simplest model that is representative of the phenomenon under consideration. 
We also remark that Eq. \eqref{eq:ode-h} has one state variable and one control variable and is similar to an RC electric circuit. The stability analysis of such RC equations is well studied in the literature and it is known that these systems are stable and controllable \cite{Webborn2017} and have good robustness properties \cite{Hao2014}. 

Utilizing the above dynamic equations, we now present three OCP models that are reflective of current practice and thermostats available. The first OCP provides a baseline that models how thermostats are operated today. The second provides the OCP for Type 1 thermostats (inherently decentralized) and the third provides the OCP for Type 2 thermostats (inherently centralized). Structurally, the OCPs are nearly the same, yielding a common algorithm discussed later.

\subsection{Baseline OCP model} \label{subsec:baseline-model}
In the baseline model, each house is equipped with a normal thermostat that implements an instantaneous control policy. At any time instant $t$, the thermostat decides to turn-on or turn-off the gas supply if the current temperature of the house is less or greater than (respectively) the set-point temperature of the thermostat. This is reflective of the thermostats that are currently deployed in the majority of households in the United States. For this baseline model, the control policy is defined by the following equation
\begin{flalign}
    q_h(t) = H\left(1-\frac{\theta_h(t)}{\bar{\bm \theta}_h}\right) \quad \forall t>0 \label{eq:baseline-control-policy}
\end{flalign}

\noindent
In Eq. \eqref{eq:baseline-control-policy}, $H(\cdot)$ is the Heaviside step function and $\bar{\bm \theta}_h$ is the set-point temperature of the thermostat. The model is a decentralized model because implementing the control policy in each household $h \in \mathcal H$ requires only local information from $h$. The baseline model is integrated over a finite time horizon $[0, T]$ to obtain a solution  $\theta^*_h(t)$ and $q^{*}_h(t)$ for every $h \in \mathcal H$. The mass flow rate  of natural gas consumed by all houses in the set $\mathcal H$ is computed from the solution as:
\begin{flalign}
d^b(t) = \sum_{h \in \mathcal H} \bm Q_h \cdot q^*_h(t)
\label{eq:total-gas-consumption-baseline} 
\end{flalign}

\noindent
In Eq. \eqref{eq:total-gas-consumption-baseline}, $d^{b}$ models the total mass flow rate of gas consumed (in \si{\kilogram\per\second}) by all houses in the set $\mathcal H$. We refer to this value as the baseline gas consumption for a set of households in $\mathcal H$.  

\subsection{Thermostat Type 1: Decentralized look-ahead OCP} \label{subsec:decentralized-lookahead} 

The decentralized look-ahead model uses ``Type 1'' smart thermostats. The output of this model is a decentralized on-off control policy for each house $h$. The control policy for each household is dependent only on current temperature and the parameters of the thermostat in that house. This model is a look-ahead model because, unlike the baseline model, this model accounts for temperature forecasts and computes a control policy accordingly. We next introduce additional notations to formulate this decentralized look-ahead OCP. 

The Type 1 thermostat has the capability to obtain a temperature forecast, $\bm \Theta(t)$, up to a set time $T$ in the future (referred to as the optimization time horizon). The first thermostat parameter is the desired temperature set-point, $\bar{\bm \theta}_h$, that the user wants to maintain inside the house $h \in \mathcal H$. The second parameter, $\bm \lambda_h \in [0, 1]$, represents the mode in which the thermostat operates. When the value of $\bm \lambda_h = 1$, the thermostat keeps the maximum deviation of the house temperature from its set point to a minimum for the entire optimization time horizon. In the other words, the DR OCP computes a control policy that minimizes the maximum deviation over the entire time horizon, defined by $\Delta_h \triangleq \max_{t\in(0, T]} |\bar{\bm \theta}_h - \theta_h(t)|$. When $\bm \lambda_h = 0$, the thermostat computes a control policy that minimizes the total natural gas consumption over the entire optimization time horizon (see Eq. \eqref{eq:gas-dec}). For all other values of $\bm \lambda_h \in (0, 1)$, the objective is to minimize a convex combination of $\Delta_h$ and $G_h$ with $\bm \lambda_h$ as the multiplier and models a trade-off between these outcomes. However, this trade-off cannot be quantified \emph{a priori}, in the sense that it is not possible to choose a value of $\bm \lambda_h$ that will guarantee a certain outcome for either $\Delta_h$ or $G_h$ unless the pareto front is solved for. 
The OCP of the decentralized look-ahead model that the smart thermostat solves for each house $h \in \mathcal H$ is formulated as follows:
\begin{subequations}
\begin{flalign}
& \min: \quad \bm \lambda_h \cdot \Delta_h + (1 - \bm \lambda_h) \cdot G_h\label{eq:obj-dec} \\ 
& \text{subject to:} \notag  \\
& \dot{\theta}_h(t) = -\bm \alpha_h \cdot (\theta_h(t) - \bm \Theta(t)) + \bm \beta_h \cdot \bm Q_h \cdot q_h(t) \notag \\ & \qquad \qquad \qquad \qquad \qquad \qquad \quad \forall t \in (0, T] \label{eq:ode-dec}, \\
& \quad G_h = \bm Q_h \cdot \int_0^T q_h(t) \label{eq:gas-dec}, \\ 
& \quad \theta_h(t=0) = \bm \theta_h^0,  \label{eq:ic-dec}\\
& \quad q_h(t) \in \{0, 1\}, \quad \forall t \in (0, T]\label{eq:control-dec} \\
& \quad \bar{\bm \theta}_h - \Delta_h  \leqslant \theta_h(t) \leqslant \bar{\bm \theta}_h + \Delta_h \quad \forall t \in (0, T]. \label{eq:deviation-dec}
\end{flalign}
\label{eq:ocp-dec}
\end{subequations}

\noindent
Here, the optimal control policy for each house is obtained by solving the OCP in \eqref{eq:ocp-dec} for each house. In the above formulation, Eq. \eqref{eq:ode-dec} represents the dynamics of heating for the house $h$, Eq. \eqref{eq:gas-dec} computes the total gas consumption for house $h$ over the entire optimization time-horizon, and Eq. \eqref{eq:ic-dec} provides the initial conditions for the dynamics in Eq. \eqref{eq:ode-dec}. Finally,  Eq. \eqref{eq:control-dec} is the binary restrictions on the control variable $q_h(t)$ and Eq. \eqref{eq:deviation-dec} provides the bounds on the state variables, $\theta_h(t)$. 

The decentralized control policy that is implemented by the thermostat allows the customer to save on natural gas consumption (incentive) without compromising on the temperature deviation from the set-point temperature. From an aggregator stand-point, it provides a reduction in peak (DR) without compromising on the comfort level of the participants. In that sense, this decentralized look-ahead model is a win-win for both the participants of the DR program and the aggregator. 
This formulation can be thought of as indirect DR, since conventional DR usually responds to prices or overall demand. Here we  have an indirect price response, in the form of reduced overall consumption. The decentralized formulation explores the limits of what can be done in the context of DR using today's smart thermostats with access to temperature and forecasts.
    
\subsection{Thermostat Type 2: Centralized look-ahead models} \label{subsec:centralized-lookahead}
We next discuss the OCP for ``Type 2'' thermostats.
Unlike a Type 1 thermostat, the Type 2 thermostat has no computational capability. The thermostat sends the temperature set-point and the current temperature of the house to the aggregator and receives and implements an on-off control policy that is computed by the aggregator. The aggregator has a forecast of the ambient temperature for the duration of the optimization time horizon and uses it to compute the control policy based on look-ahead information. The model is centralized because the aggregator utilizes (i) set-point temperatures of each household, (ii) initial temperatures of each household, (iii) parameters of each house, (iv) ambient temperature forecasts, and (v) additional constraints on consumption to solve a single OCP that computes a control policy for each household. The computed control policy is transmitted to the thermostat in each household for implementation. 

Similar to the decentralized look-ahead models, there is an optimization time horizon $[0, T]$. We let, $\Delta_h$, denote the maximum deviation between the set-point temperature and the actual temperature of house $h$ over the entire optimization time horizon. To model peak shaving, we let $\bm D$ denote the maximum total amount of gas consumed by all the households in $\mathcal H$ and let $\bm \gamma \in (0, 1]$ denote the percentage of gas reduction that the aggregator achieves with the centralized DR program.  With these notations, we present the following look-ahead model with two different objective functions (i) mean deviation objective and (ii) max deviation objective:
\begin{subequations}
\begin{flalign}
\text{Mean Deviation: } & \min \quad \frac{1}{|\mathcal H|}  \sum_{h \in \mathcal H} \Delta_h \label{eq:mean} \\ 
\text{Max. Deviation: } & \min \quad \left(\max_{h \in \mathcal H} \Delta_h\right) \label{eq:max}
\end{flalign}
\label{eq:objectives}
\end{subequations}

\noindent
In both the objective functions, the value of $\Delta_h$ is interpreted as discomfort of participants. The objective function in Eq. \eqref{eq:mean} and Eq. \eqref{eq:max}  minimizes the average and maximum discomfort, respectively, over all the households participating in the DR program. The centralized look-ahead OCP model with objective functions in Eq. \eqref{eq:objectives} is subject to the following constraints:
\begin{subequations}
\begin{flalign}
    & \quad \dot{\theta}_h(t) = -\bm \alpha_h \cdot (\theta_h(t) - \bm \Theta(t)) + \bm \beta_h \cdot \bm Q_h \cdot q_h(t) \notag \\ 
    & \qquad \qquad \qquad \qquad \qquad \qquad \quad \forall t \in (0, T], \forall h \in \mathcal H  \label{eq:ode-cen}, \\
    & \quad \theta_h(t=0) = \bm \theta_h^0,  \quad \forall h \in \mathcal H \label{eq:ic-cen}\\
    & \quad q_h(t) \in \{0, 1\}, \quad \forall t \in (0, T], \forall h \in \mathcal H \label{eq:control-cen} \\
    & \quad \bar{\bm \theta}_h - \Delta_h  \leqslant \theta_h(t) \leqslant \bar{\bm \theta}_h + \Delta_h \notag \\ 
    & \qquad \qquad \qquad \qquad \qquad \qquad \quad \forall t \in (0, T], \forall h \in \mathcal H, \label{eq:deviation-cen} \\
    & \quad \sum_{h \in \mathcal H} \bm Q_h \cdot q_h(t) \leqslant \bm \gamma \cdot \bm D \quad \forall t \in (0, T]. \label{eq:peak-shaving}
\end{flalign}
\label{eq:centralized-constraints}
\end{subequations}

\noindent
The constraints in Eq. \eqref{eq:peak-shaving} enable the operator to enforce a peak mass flow rate of $\bm \gamma \cdot \bm D$. Here, $\bm D$ is 
determined by the system operator by estimating of how much total gas is consumed by all the households during winter peaks. $\bm \gamma$ is a percentage factor that the operator can vary to perform reliability analysis in the case of extreme events. We refer to $\bm \gamma$ as the curtailment factor that indicates the percentage by which the operator wishes to curtail the maximum mass flow rate in the natural gas system. The centralized model offers an advantage over the decentralized one in that once $\bm \gamma$ is chosen, the operator already knows what to expect regarding peak consumption, and if the deviation performance corresponding to a chosen value of $\bm \gamma$ is unsatisfactory, then it is clear that a lower value of $\bm \gamma$ needs to be considered to achieve better deviation performance. 

%
%


\section{Algorithms} \label{sec:algo}

This section presents a receding-horizon algorithm to solve the OCPs described in section \ref{sec:overview}. To solve the OCPs, we first discretize the state and control variables over time. It is useful to separate the temporal resolution of the control variables from the state variables to better model the capabilities of the controller and limit chattering
\cite{slotine1991applied}.
Before we present the receding-horizon algorithm, we first present a discretized version of the decentralized and centralized look-ahead models. 

In the algorithm, the optimization time horizon $[0, T]$ is discretized into $n$ equispaced time instants $\mathcal K = \{0=t_0, t_1 =\bm \Delta \bm t, \dots, t_n = n\bm \Delta \bm t = T\}$ with a time step of $\bm \Delta \bm t = T/ n$, for every house $h \in \mathcal H$. For OCP's based on Type 1 thermostats, the decentralized look-ahead model is a Mixed-Integer Linear Program (MILP) of the form:

\begin{subequations}
\begin{flalign}
    & \min: \quad \bm \lambda_h \cdot \Delta_h + (1 - \bm \lambda_h) \cdot G_h \text{\quad subject to:} \label{eq:obj-milp1} \\ 
    & \quad \frac{\theta_h(t_{k+1}) - \theta_h(t_{k})}{\bm \Delta \bm t} =  -\bm \alpha_h (\theta_h(t_k) - \bm \Theta(t_k)) + 
    \bm \beta_h \cdot \bm Q_h \cdot q_h(t_k), \; \forall k \in \{0, \cdots, n-1\}  \label{eq:ode-milp1} \\
& \quad \theta_h(t_0) = \bm \theta_h^0, \label{eq:ic-milp1} \\
& \quad G_h =  \bm Q_h \cdot  \sum_{k=0}^n q_h(t_k), \label{eq:heat-flux-gas-milp1}  \\
& \quad \bar{\bm \theta}_h - \Delta_h  \leqslant \theta_h(t_k) \leqslant \bar{\bm \theta}_h + \Delta_h \quad \forall k \in \mathcal K, \label{eq:deviation-milp1} \\
& \quad q_h(t_k) \in \{0, 1\} \quad \forall k \in \mathcal K \label{eq:binary-milp1}
\end{flalign}
\label{eq:milp-1}
\end{subequations}

The MILPs for the Type 2 OCP are solved as a single optimization problem for all houses by the aggregator. The constraints of these MILPs after discretization, are as follows:

\begin{subequations}
\begin{flalign}
& \frac{\theta_h(t_{k+1}) - \theta_h(t_{k})}{\bm \Delta \bm t} =  -\bm \alpha_h (\theta_h(t_k) - \bm \Theta(t_k)) + \bm \beta_h \cdot \bm Q_h \cdot q_h(t_k), \; \forall k \in \{0, \cdots, n-1\} \; \forall h \in \mathcal H \label{eq:ode-milp2} \\
& \theta_h(t_0) = \bm \theta_h^0, \quad \forall h \in \mathcal H \label{eq:ic-milp2} \\
& \bar{\bm \theta}_h - \Delta_h  \leqslant \theta_h(t_k) \leqslant \bar{\bm \theta}_h + \Delta_h \quad \forall k \in \mathcal K, \quad \forall h \in \mathcal H, \label{eq:deviation-milp2}\\
& \sum_{h \in \mathcal H} \bm Q_h \cdot q_h(t_k) \leqslant \bm \gamma \cdot \bm D \quad \forall k \in \mathcal K, \text{ and } \label{eq:peak_reduction-milp2} \\
& q_h(t_k) \in \{0, 1\} \quad \forall k \in \mathcal K, \quad \forall h \in \mathcal H. \label{eq:binary-milp2}
\end{flalign}
\label{eq:constraints-milp2}
\end{subequations}


\subsection{A Receding-Horizon (RH) algorithm} \label{subsec:rh}
\begin{figure*}[htb]
    \centering
    \includegraphics[scale=0.2]{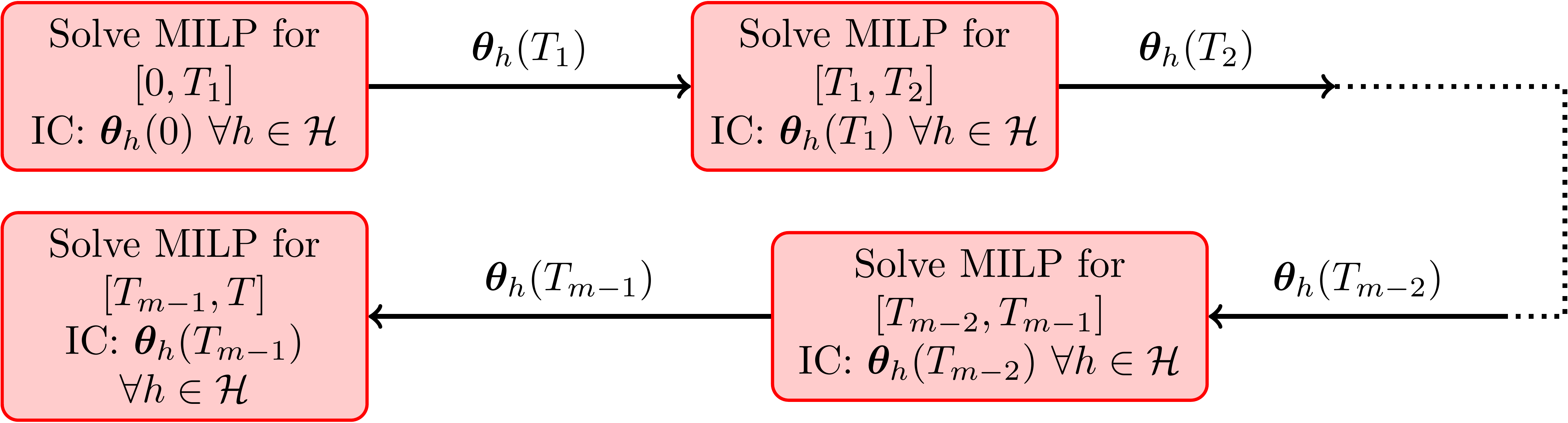}
    \caption{Schematic of the RH algorithm. IC denotes initial conditions.}
    \label{fig:rh-schematic}
\end{figure*}

One issue that hinders solving the OCP MILPs to global optimality is the size of the MILPs. For the Type 1 thermostat there is one MILP for each house, $h \in \mathcal H$, that is solved separately. The size of the MILP (Eqs. \eqref{eq:milp-1} and \eqref{eq:constraints-milp2}) is large when $\bm \Delta \bm t$ is small  or the optimization time horizon $T$ is large. To address the computational issue associated with large MILPs, we use a Receding Horizon (RH) algorithm. The approach decomposes the time horizon into smaller time horizons and solves each new problem sequentially. More formally, the RH splits the time horizon $[0, T]$ into $[0, T_1]$, $[T_1, T_2]$, $\dots$, $[T_{m-1}, T_m]$ where $T_m = T$. Each of the resulting $m$ MILPs for each smaller time horizon is solved sequentially. Here, the initial condition for the $i$\textsuperscript{th} MILP is determined by the solution of the $(i-1)$\textsuperscript{th} MILP at the final discrete time instant. The schematic of the RH algorithm is shown in Fig. \ref{fig:rh-schematic}. 

\section{Case Studies} \label{sec:case-study}
In this section, we present two case studies. The first case study analyzes the efficacy of the models for a typical winter day in the Chicago area and the second case study analyzes the impact of the 2019 Polar vortex  in the Chicago area from temperature data. For both the case studies, the set of houses and the parameters for each house remain unchanged. The Julia Programming language \cite{bezanson2017julia} was used for all the implementations and Gurobi \cite{bixby2007gurobi} was used to solve the MILPs that were obtained by discretizing the OCPs. All the computational experiments were run on an Intel Broadwell E5-2695 processor with a base clock rate of 2.10 GHz and a RAM of 128 GB.

\subsection{House data} \label{subsec:house-data}
For all the case studies, the total number of houses in the set $\mathcal H$ is set to $140$. $14$ categories of houses were chosen ranging from single bedroom apartments to single family homes. The thermal coefficients for each house was calculated according to the formulae derived in Appendix \ref{sec:thermal-coeff}. The furnace burn rate, $\bm Q_h$, for each house was set in the range of $[3\times 10^{-5}, 12\times 10^{-5}]$ \si{\kilogram\per\second} depending on the square footage of the corresponding house. The value of the parameters specific to each house category (volume, square footage, initial temperature, set-point temperature, etc.) are made available at \url{https://github.com/kaarthiksundar/ng-dr-data}. The other parameter values used for computing the thermal coefficients for each house are shown in Table \ref{tab:parameters}. 
\begin{table}[htb]
    \centering
    \begin{tabular}{c|c}
        \toprule 
        Parameter & Value (SI Units) \\
        \midrule
        $\bm C^a$ & $0.718 \times 10^3$ (\si{\joule\per\kilogram\per\kelvin}) \\
        $\bm \rho^a$ & $1.2754$ (\si{\kilogram\per\cubic\meter}) \\
        $\bm \kappa_h$ $\forall h \in \mathcal H$ & $0.11$ (\si{\watt\per\square\meter\per\kelvin}) \\
        $\bm E^g$ & $45938 \times 10^3$ (\si{\joule\per\kilogram}) \\
        \bottomrule
    \end{tabular}
    \caption{Parameter values used for the case studies.}
    \label{tab:parameters}
\end{table}

\subsection{Ambient temperature data} 
\begin{figure}[htbp]
\centering
    \begin{subfigure}[t]{0.4\textwidth}
    \centering
    \includegraphics[scale=0.5]{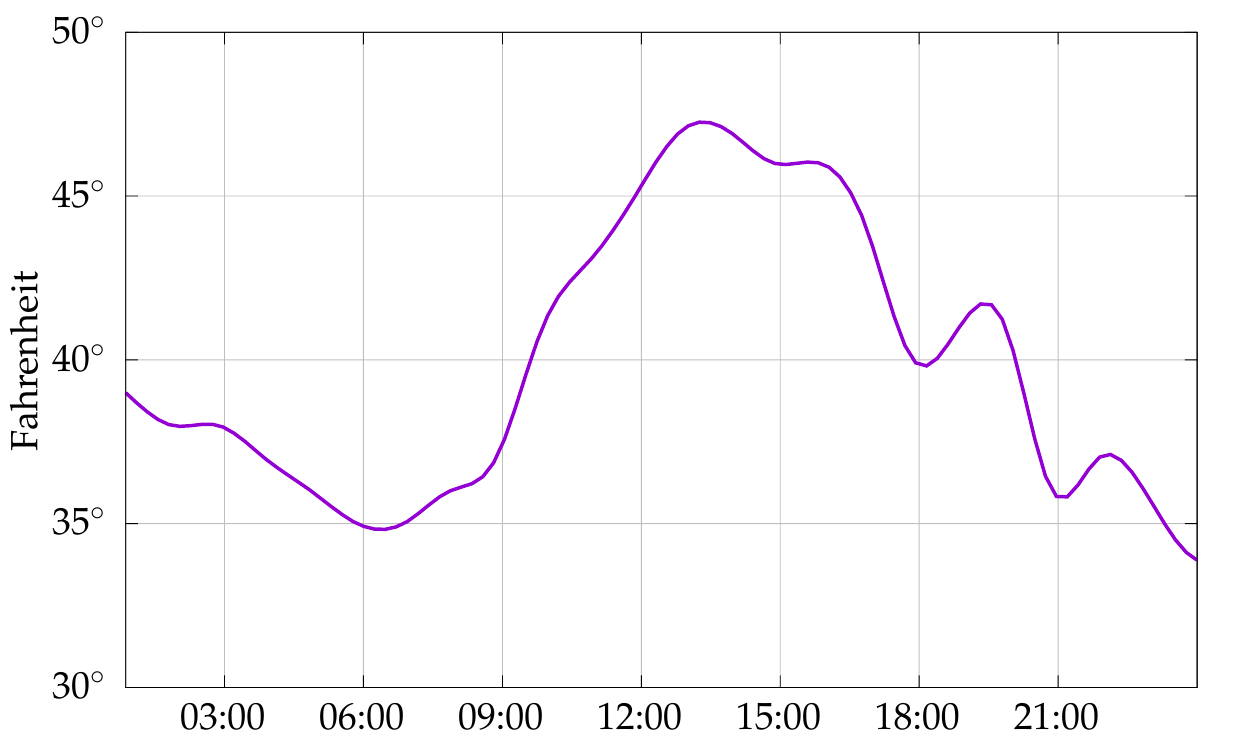}
    \caption{Ambient temperature data for the Chicago area on Jan. 2, 2020. This forecast is representative of a typical winter day in the Chicago area.}
    \label{fig:ambient_typical}
    \end{subfigure}
    \hspace{3ex}
    \begin{subfigure}[t]{0.4\textwidth}
        \centering
        \includegraphics[scale=0.5]{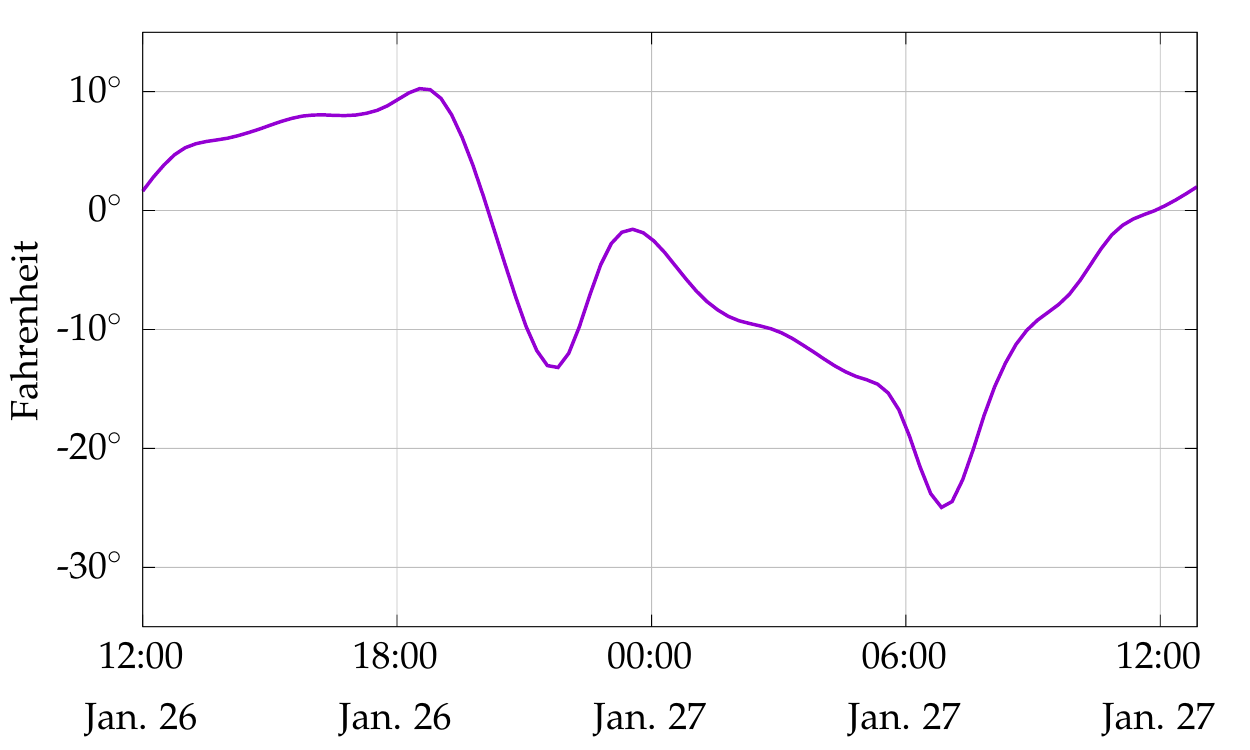}
    \caption{Ambient temperature data for the Chicago area on Jan. 26-27, 2019. This is the day when the effect of the polar vortex was most severe.}
    \label{fig:ambient_pv}
    \end{subfigure}
    \caption{Ambient temperature data for both the case studies in this article. Note that the temperature varies by at most 15 degrees for the typical day, while the range is almost 40 degrees for the  day of the polar vortex.}
    \label{fig:ambient}
\end{figure}

For a typical winter day case study, historical ambient temperature data for the Chicago area on January 2, 2020 was chosen. For the polar vortex case, the temperature data from January 26, 2019 noon to January 27, 2019 noon was chosen. The temperature data for both cases were obtained from \url{https://www.noaa.gov/weather}. The plot of the temperature profiles for both the cases are shown in Figures \ref{fig:ambient_typical} and \ref{fig:ambient_pv}, respectively. Throughout the rest of the article, we refer to the two case studies as ``typical-day'' and ``polar-vortex''.

\subsection{Results from the Type 1 thermostat model} \label{subsec:dec-results}

Here we present the results obtained for the the Type 1 thermostat model in Eq. \eqref{eq:ocp-dec} for both case studies and compare it with the results from the baseline model in Eq. \eqref{eq:dynamics} and Eq. \eqref{eq:baseline-control-policy}. For the decentralized look-ahead model, the value of $\bm \lambda_h$ for all the houses is varied in the set $\{0.65, 0.7, 0.75, 0.8, 0.85, 0.9, 0.95, 1\}$ and the look-ahead time period, denoted by $T_{rh}$ (in hours) is chosen from the set $\{1, 3\}$. The look-ahead OCP for each house is solved separately but in parallel. The state and control discretization for the MILP reformulations for the look-ahead OCPs were set to $1$ and $3$ minutes, respectively.

\begin{table}[htbp]
    \centering
    \begin{tabular}{ccccc}
        \toprule
        \multirow{3}{*}{$\bm \lambda_h$} & \multicolumn{2}{c}{Computation Time (seconds)} & \multicolumn{2}{c}{Computation Time (seconds)} \\
        & \multicolumn{2}{c}{``typical-day'' case} & \multicolumn{2}{c}{``polar-vortex'' case}\\
        \cmidrule{2-5}
        & $T_{rh}=1$ hr & $T_{rh}=3$ hr & $T_{rh}=1$ hr & $T_{rh}=3$ hr \\ 
        \midrule 
        \begin{filecontents*}{tables/dec-times.csv}
lamba,rhonet,rhthreet,rhonep,rhthreep
0.65,0.06,4.76,0.07,19.21
0.70,0.05,2.75,0.07,19.52
0.75,0.05,1.78,0.07,18.41
0.80,0.05,1.45,0.06,24.38
0.85,0.05,1.23,0.06,20.49
0.90,0.04,1.09,0.05,20.08
0.95,0.04,0.93,0.05,20.57
1.00,0.03,0.20,0.03,0.19
\end{filecontents*}
\csvreader[late after line=\\]{tables/dec-times.csv}{1=\k,2=\rhonet,3=\rhthreet,4=\rhonep,5=\rhthreep}{\k & \rhonet & \rhthreet & \rhonep & \rhthreep}
        \bottomrule
    \end{tabular}
    \caption{Computation time in seconds for the decentralized look-ahead model averaged over all the houses. As $T_{rh}$ increases, it is evident that the size of the MILP increases and hence, an increase in run time is observed for both the case studies as we move from $T_{rh} = 1$ hour to $T_{rh} = 3$ hours. In the text, we have explained why runs were not performed for a lower value of $\bm \lambda_h$.}
    \label{tab:dec-times}
\end{table}

Table \ref{tab:dec-times} shows the average computation time, in seconds, to solve the MILPs to optimality. The average is computed over all the houses since each house has its own MILP. As $T_{rh}$ increases, it is evident that the size of the MILP increases and run time increases. 

\begin{table}[h!]
\begin{subtable}{\textwidth}
    \centering
    \begin{tabular}{ccrrrrcc}
        \toprule 
        \multirow{2}{*}{model} & \multirow{2}{*}{$\bm \lambda_h$} & \multicolumn{4}{c}{deviation (\si{\fahrenheit})} & \multirow{2}{3cm}{\centering gas consumption (\si{\kilogram})} & \multirow{2}{3cm}{\centering consumption per house (\si{\kilogram})}\\
        \cmidrule{3-6}
        & & mean & max. & min. & std. dev. & & \\
        \midrule 
        \begin{filecontents*}{tables/dec-dev-typ-1.csv}
model,lambda,mean,max,min,std,gas,gasperhouse
baseline,NA,3.28,4.83,1.37,0.91,235.1538,1.6797
decentralized,0.65,7.71,36.99,1.38,10.92,155.5983,1.1114
decentralized,0.70,3.45,9.82,1.19,1.80,217.5768,1.5541
decentralized,0.75,2.74,3.60,1.19,0.65,225.3636,1.6097
decentralized,0.80,2.51,3.60,1.06,0.58,226.7347,1.6195
decentralized,0.85,2.38,3.40,1.06,0.55,227.2455,1.6232
decentralized,0.90,2.25,3.33,1.05,0.54,227.9151,1.6280
decentralized,0.95,2.17,3.33,1.05,0.54,228.2083,1.6301
decentralized,1.00,2.14,3.33,1.02,0.55,228.3865,1.6313
\end{filecontents*}
\csvreader[late after line=\\]{tables/dec-dev-typ-1.csv}{1=\model,2=\l,3=\mean,4=\max,5=\min,6=\stddev,7=\gas,8=\perhouse}{\model & \l & \mean & \max & \min & \stddev & \gas & \perhouse}
        \bottomrule
    \end{tabular}
    \caption{``Typical-day'' case study with $T_{rh}=1$ hour}
    \label{tab:rh1_typ}
\end{subtable} \\\vspace{3ex}

\begin{subtable}{\textwidth}
\centering 
\begin{tabular}{ccrrrrcc}
        \toprule 
        \multirow{2}{*}{model} & \multirow{2}{*}{$\bm{\lambda}_h$} & \multicolumn{4}{c}{deviation (\si{\fahrenheit})} & \multirow{2}{3cm}{\centering gas consumption (\si{\kilogram})} & \multirow{2}{3cm}{\centering consumption per house (\si{\kilogram})}\\
        \cmidrule{3-6}
        & & mean & max. & min. & std. dev. & & \\
        \midrule 
        \begin{filecontents*}{tables/dec-dev-typ-3.csv}
model,lambda,mean,max,min,std,gas,gasperhouse
baseline,NA,3.28,4.83,1.37,0.91,235.1538,1.6797
decentralized,0.65,2.38,3.48,1.04,0.56,226.5980,1.6186
decentralized,0.70,2.31,3.48,1.04,0.55,227.2082,1.6229
decentralized,0.75,2.26,3.41,1.04,0.53,227.5393,1.6253
decentralized,0.80,2.21,3.33,1.04,0.54,227.8233,1.6273
decentralized,0.85,2.17,3.33,1.04,0.54,227.9923,1.6285
decentralized,0.90,2.15,3.33,1.04,0.55,228.2191,1.6301
decentralized,0.95,2.14,3.33,1.02,0.55,228.3163,1.6308
decentralized,1.00,2.14,3.33,1.02,0.55,228.3811,1.6313
\end{filecontents*}
\csvreader[late after line=\\]{tables/dec-dev-typ-3.csv}{1=\model,2=\l,3=\mean,4=\max,5=\min,6=\stddev,7=\gas,8=\perhouse}{\model & \l & \mean & \max & \min & \stddev & \gas & \perhouse}
        \bottomrule
    \end{tabular}
    \caption{``Typical-day'' case study with $T_{rh}=3$ hour}
    \label{tab:rh3_typ}
\end{subtable}
\caption{Statistics of temperature deviation and amount of gas consumed for the entire day for the ``typical-day'' case study. Note that reduction in gas consumption does not come at the cost of compromising on the deviation performance (in comparison to the baseline model). Thus, utilizing smart thermostats with look-ahead features can lead to substantial amount of reduction in consumption  without compromising the deviation performance for a typical winter day.}
\label{tab:typ}
\end{table}

\begin{table}[h!]
    \begin{subtable}{\textwidth}
    \centering
    \begin{tabular}{ccrrrrcc}
        \toprule 
        \multirow{2}{*}{model} & \multirow{2}{*}{$\bm \lambda_h$} & \multicolumn{4}{c}{deviation (\si{\fahrenheit})} & \multirow{2}{3cm}{\centering gas consumption (\si{\kilogram})} & \multirow{2}{3cm}{\centering consumption per house (\si{\kilogram})}\\
        \cmidrule{3-6}
        & & mean & max. & min. & std. dev. & & \\
        \midrule 
        \begin{filecontents*}{tables/dec-dev-pv-1.csv}
model,lambda,mean,max,min,std,gas,gasperhouse
baseline,NA,5.97,14.70,2.41,3.70,560.9169,4.0065
decentralized,0.65,18.03,89.76,1.86,28.25,398.2451,2.8446
decentralized,0.70,7.03,20.21,1.74,4.88,545.9659,3.8998
decentralized,0.75,5.79,14.70,1.68,3.88,558.9702,3.9926
decentralized,0.80,5.53,14.70,1.63,3.99,562.1605,4.0154
decentralized,0.85,5.38,14.70,1.53,4.08,563.3577,4.0240
decentralized,0.90,5.30,14.70,1.53,4.14,563.9166,4.0280
decentralized,0.95,5.25,14.70,1.53,4.18,564.4836,4.0320
decentralized,1.00,5.23,14.70,1.53,4.19,566.9800,4.0499
\end{filecontents*}
\csvreader[late after line=\\]{tables/dec-dev-pv-1.csv}{1=\model,2=\l,3=\mean,4=\max,5=\min,6=\stddev,7=\gas,8=\perhouse}{\model & \l & \mean & \max & \min & \stddev & \gas & \perhouse}
        \bottomrule
    \end{tabular}
    \caption{``Polar-vortex'' case study with $T_{rh}=1$ hour}
    \label{tab:rh1_pv}
    \end{subtable} \\\vspace{3ex}
    
    \begin{subtable}{\textwidth}
    \centering
    \begin{tabular}{ccrrrrcc}
        \toprule 
        \multirow{2}{*}{model} & \multirow{2}{*}{$\bm \lambda_h$} & \multicolumn{4}{c}{deviation (\si{\fahrenheit})} & \multirow{2}{3cm}{\centering gas consumption (\si{\kilogram})} & \multirow{2}{3cm}{\centering consumption per house (\si{\kilogram})}\\
        \cmidrule{3-6}
        & & mean & max. & min. & std. dev. & & \\
        \midrule 
        \begin{filecontents*}{tables/dec-dev-pv-3.csv}
model,lambda,mean,max,min,std,gas,gasperhouse
baseline,NA,5.97,14.70,2.41,3.70,560.9169,4.0065
decentralized,0.65,5.51,14.82,1.60,4.18,559.5426,3.9967
decentralized,0.70,5.43,14.82,1.53,4.18,560.2889,4.0021
decentralized,0.75,5.36,14.75,1.53,4.20,560.9531,4.0068
decentralized,0.80,5.31,14.70,1.53,4.21,561.5633,4.0112
decentralized,0.85,5.25,14.70,1.53,4.18,562.8647,4.0205
decentralized,0.90,5.24,14.70,1.53,4.19,562.9878,4.0213
decentralized,0.95,5.23,14.70,1.53,4.19,563.2481,4.0232
decentralized,1.00,5.22,14.70,1.53,4.19,564.7536,4.0340
\end{filecontents*}
\csvreader[late after line=\\]{tables/dec-dev-pv-3.csv}{1=\model,2=\l,3=\mean,4=\max,5=\min,6=\stddev,7=\gas,8=\perhouse}{\model & \l & \mean & \max & \min & \stddev & \gas & \perhouse}
        \bottomrule
    \end{tabular}
    \caption{``Polar-vortex'' case study with $T_{rh}=3$ hour}
    \label{tab:rh3_pv}
    \end{subtable}
    \caption{Statistics of temperature deviation and amount gas of consumed for the entire day for the ``polar-vortex'' case study.  We observe that though the reduction in gas consumption  is not as substantial as in the ``typical-day'' case study, there exist regimes ($\bm \lambda_h$ values) for which the decentralized look-ahead models provide some reduction in gas consumption along with good deviation performance.}
\label{tab:pv}
\end{table}

Tables \ref{tab:typ} and \ref{tab:pv} show the statistics of temperature deviation and the total gas consumption obtained using the Type-1 thermostat for $T_{rh} \in \{1, 3\}$ hours. The deviation performance is a measure of discomfort caused to the households, i.e., the greater the deviation, the greater the discomfort. The level of discomfort for $\bm \lambda_h = 0.65$ provides ample evidence that experiments for $\bm \lambda_h < 0.65$ are unnecessary.
In all the tables, the first row is computed by solving the baseline model. The results indicate that a lower value of the look-ahead time period results in more short-sightedness of the solutions. Reducing the value of $T_{rh}$ degrades the deviation performance needed to achieve lower gas consumption. This trend is observed in both the Tables \ref{tab:typ} and \ref{tab:pv}. 

For the ``typical-day'' case study results in Tables \ref{tab:rh1_typ} and \ref{tab:rh3_typ}, we observe reduction in gas consumption per house for many values of $\bm \lambda_h$. The values of the gas consumption per house indirectly show more aggregation benefits as well as economic benefits since savings in gas consumption directly results in cost savings for the consumer. For instance, in Table \ref{tab:rh1_typ}, when the value of $\bm \lambda_h = 0.75$, the gas savings per house is $0.07$ \si{\kilogram}. When the number of houses aggregated is around $10000$, this leads to a total gas consumption reduction of $\approx 700$ \si{\kilogram}, which is a substantial amount of reduction in consumption even on a typical winter day. Furthermore, it is important to note that this reduction in gas consumption does not come at the cost of compromising on the deviation performance (in comparison to the baseline model).  In summary, utilizing smart thermostats with look-ahead features reduces a substantial amount of consumption without compromising deviation performance for a typical winter day. 

In the ``polar-vortex'' case study, we observe that though reduction in gas consumption  is not as substantial as in the ``typical-day'' case study, there exist regimes ($\bm \lambda_h$ values) for which the decentralized look-ahead models provide  some reduction in gas consumption along with good deviation performance. Also, unlike the ``typical-day'' case, where consumption reduction is achieved for all $\bm \lambda_h$ values that provide deviation performance comparable to the baseline model, the values of $\bm \lambda_h$ closer to  $1$ cause an increase in the gas consumption. This is again expected due to the severity of the polar vortex. Nevertheless, if the operator can ensure all the smart thermostats for the decentralized model have appropriate values of $\bm \lambda_h$, it can still lead to some reduction in gas consumption.

Finally, the Tables \ref{tab:pv-gas-savings} and \ref{tab:typ-gas-savings} present the amount of reduction in gas consumption that can be achieved for both the ``polar-vortex'' and the ``typical-day'' case studies, respectively, for certain regimes ($\bm \lambda_h$ values) when compared against the baseline model which is reflective of current practice. The value of $T_{rh}$ was set to $3$ hours for the results in both the tables. The results show that even in an extreme event scenario like the polar vortex, there exists certain regimes of operation of the thermostats that provide an overall reduction in natural gas consumption while providing a better temperature deviation performance than what can be achieved by current practice (baseline model). In a typical winter day, this reduction is amplified even further with a wider range of $\bm \lambda_h$ values providing the same type of benefit to the customers. From an aggregator stand-point, installation of smart Type 1 thermostats on a greater number of houses only increases the aggregation benefits as shown in the final column of the Tables in \ref{tab:savings}.
\begin{table}[htbp]
    \begin{subtable}{\textwidth}
    \centering
    \begin{tabular}{ccc}
        \toprule 
        $\bm \lambda_h$ & reduction (\si{\kilogram}) & reduction for 10000 houses (\si{\kilogram})\\
        \midrule
        \begin{filecontents*}{tables/dec-savings-pv.csv}
lambda,savings,savingsfortenthousand
0.65,1.3743,98.16
0.70,0.6280,44.85
\end{filecontents*}
\csvreader[late after line=\\]{tables/dec-savings-pv.csv}{1=\l,2=\gas,3=\gasten}{\l & \gas & \gasten}
        \bottomrule
    \end{tabular}
    \caption{``Polar-vortex'' case study}
    \label{tab:pv-gas-savings}
    \end{subtable} \\\vspace{3ex}
    
    \begin{subtable}{\textwidth}
    \centering
    \begin{tabular}{cccc}
        \toprule 
        $\bm \lambda_h$ & reduction (\si{\kilogram}) & reduction for 10000 houses (\si{\kilogram})\\
        \midrule
        \begin{filecontents*}{tables/dec-savings-typ.csv}
lambda,savings,savingsfortenthousand
0.65,8.5558,611.13
0.70,7.9456,567.54
0.75,7.6145,543.89
0.80,7.3305,523.61
0.85,7.1615,511.54
0.90,6.9347,495.34
0.95,6.8375,488.39
1.00,6.7727,483.76
\end{filecontents*}
\csvreader[late after line=\\]{tables/dec-savings-typ.csv}{1=\l,2=\gas,3=\gasten}{\l & \gas & \gasten}
        \bottomrule
    \end{tabular}
    \caption{``Typical-day'' case study}
    \label{tab:typ-gas-savings}
    \end{subtable} 
    \caption{Regimes (values of $\bm \lambda_h$) for which the decentralized model provides considerable reduction in gas consumption in comparison for the baseline model for both the ``polar-vortex'' and the ``typical-day'' case studies. The aggregation benefits are shown in the last column of both the tables when this model is used in a locality of 10000 houses. }
    \label{tab:savings}
\end{table}

\subsection{Results for the Type-2 Thermostat Models} \label{subsec:cen-results}

Given the Type 1 thermostat model results, we focus on Type 2 thermostat results for the polar-vortex case study.
For these results, the following parameters are used: (i) $\bm D$ (peak mass flow rate) is set to the maximum mass flow rate of the baseline model and (ii) $\bm \gamma$ (curtailment factor) is chosen from the set $\{0.85, 0.90, 0.95\}$. Similar to the previous results, all the experiments for the model are run for $T_{rh} \in \{1, 3\}$ hours. In all the results, ``min.'' is used to denote the Type-2 thermostat model with the mean deviation objective and ``min. max.'' is used to denote the max. deviation objective. Finally, a time limit of 1.5 hours was placed on every computational experiment.

\begin{table}[htbp]
    \centering
    \begin{tabular}{ccccc}
    \toprule
        \multirow{2}{*}{$\bm \gamma$} & \multicolumn{2}{c}{Avg. computation time (sec)} & \multicolumn{2}{c}{Avg. optimality gap (\%)} \\
    \cmidrule{2-5} 
    & $T_{rh} = 1$ & $T_{rh} = 3$ & $T_{rh} = 1$ & $T_{rh} = 3$\\ 
    \midrule 
    \multicolumn{5}{l}{min. model} \\
    \begin{filecontents*}{tables/min.csv}
gamma,rhone,rhthree,rhone,rhthree
0.85,TO,TO,8.87,26.50
0.90,TO,TO,6.55,22.70
0.95,TO,TO,5.98,19.54
\end{filecontents*}
\csvreader[late after line=\\]{tables/min.csv}{1=\gamma,2=\tone,3=\tthree,4=\oone,5=\othree}{\gamma & \tone & \tthree & \oone & \othree}
    \midrule
    \multicolumn{5}{l}{min. max. model} \\
    \begin{filecontents*}{tables/min_max.csv}
gamma,rhone,rhthree,rhone,rhthree
0.85,14.94,1590.21,0,0
0.90,14.19,948.91,0,0
0.95,13.31,766.85,0,0
\end{filecontents*}
\csvreader[late after line=\\]{tables/min_max.csv}{1=\gamma,2=\tone,3=\tthree,4=\oone,5=\othree}{\gamma & \tone & \tthree & \oone & \othree}
    \bottomrule
    \end{tabular}
    \caption{Average computation times and optimality gaps for Type-2 thermostat models. Here, ``TO'' indicates that computation did not complete within the allotted time limit of 1.5 hours.  The high computation time and optimality gap values for the ``min.'' model when compared against the ``min. max.'' model are explained by degeneracy in the solution space of the OCPs, i.e., existence of multiple solutions with the same mean deviation when minimizing mean deviation.
    This degeneracy is removed  when the OCP minimizes the maximum deviation and the computation time improves. Also, as $T_{rh}$ increases, the number of variables and constraints also increases and naturally increases the computation times.}
    \label{tab:cen-times}
\end{table}

Table \ref{tab:cen-times} shows the average computation time and the average optimality gap for all the runs of the Type 2 thermostat model with the two different objectives. In Table \ref{tab:cen-times}, the averages are computed over each run of the centralized MILP in Sec. \ref{subsec:rh}.   One immediate observation is the high computation time and optimality gap values for the ``min.'' model when compared against the ``min. max.'' model. This trend can be explained by degeneracy in the solution space of the OCPs, i.e., existence of multiple solutions with the same mean deviation when minimizing mean deviation.
This degeneracy is removed when minimizing the maximum deviation which reduces the computation time. Also, as $T_{rh}$ increases, the number of variables and constraints increases which leads to increased computation times. Finally, as the value of the curtailment factor ${\bm \gamma}$ decreases, the control policy bounds natural gas consumption. In an extreme event setting, this can make finding optimal solutions more difficult. Physically, this is equivalent to constraining the capacity of the system. 

\begin{figure}[htbp]
\begin{subfigure}{\textwidth}
    \centering
    \includegraphics[scale=0.8]{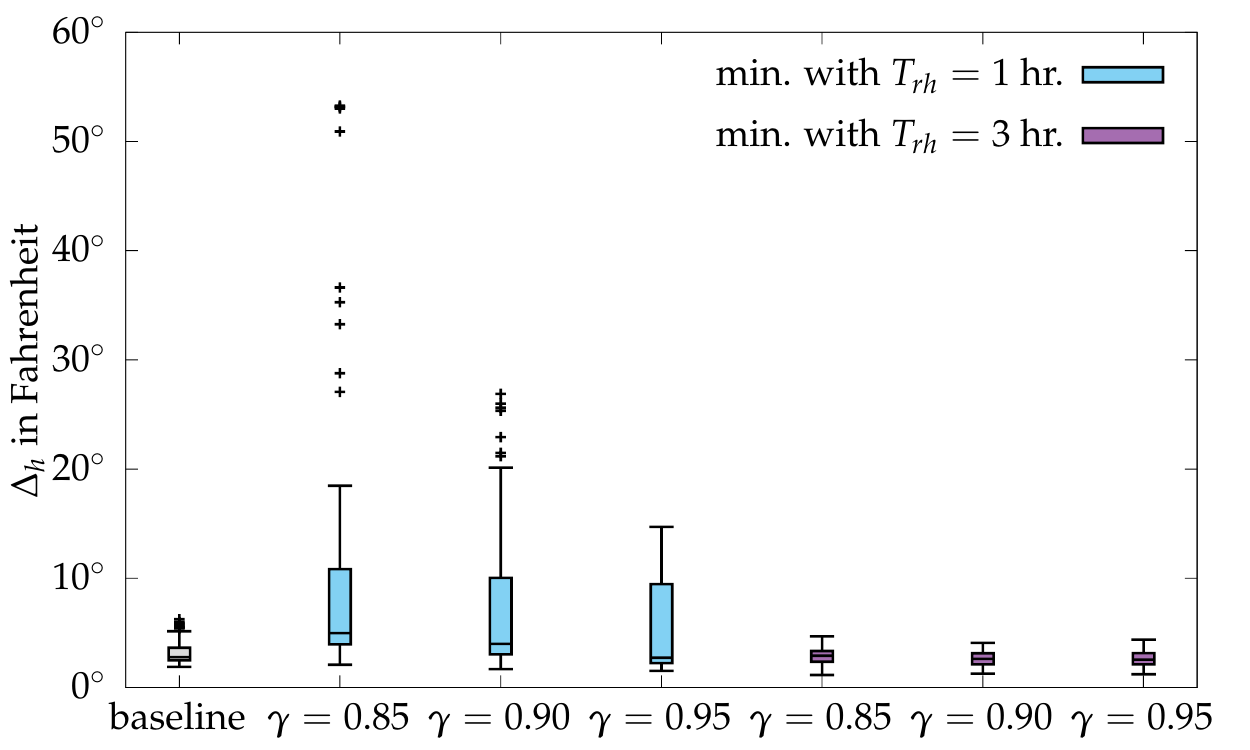}
    \caption{Centralized look-ahead model with min. objective function.}
    \label{fig:min-box}
\end{subfigure} \\\vspace{3ex} 

\begin{subfigure}{\textwidth}
    \centering
    \includegraphics[scale=0.8]{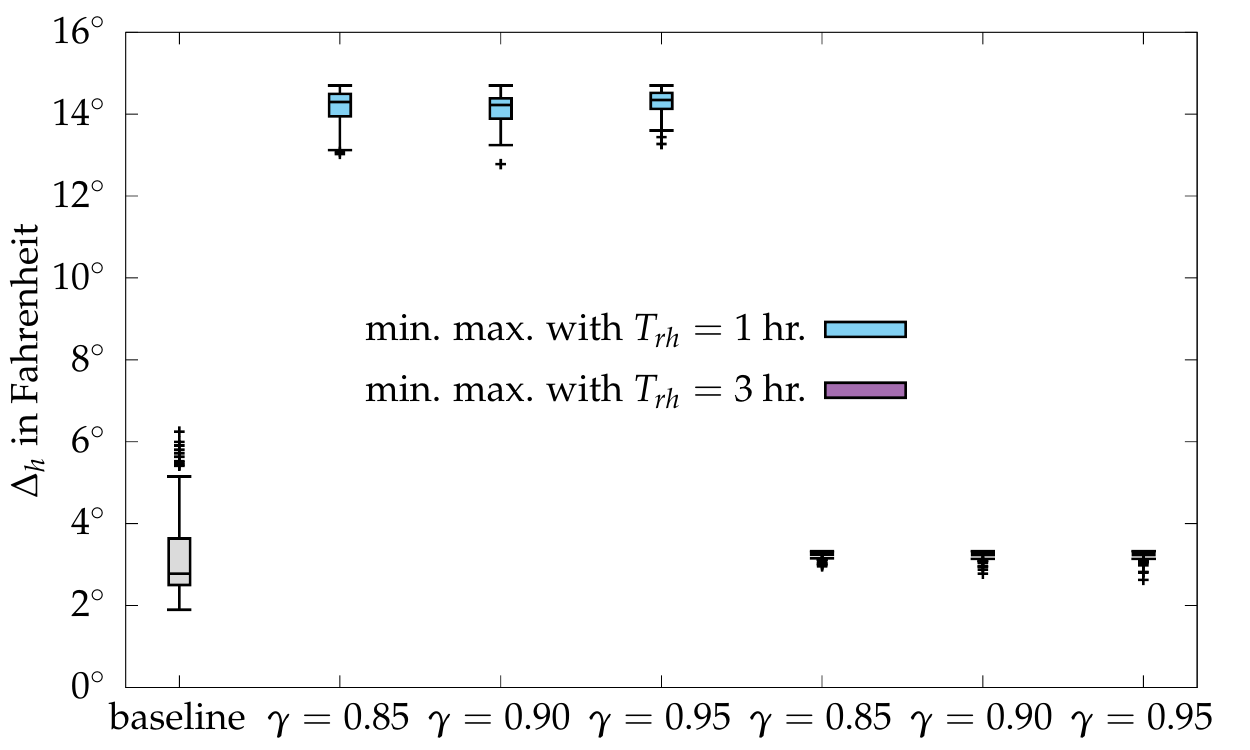}
    \caption{Centralized look-ahead model with min. max. objective function}
    \label{fig:minmax-box}
\end{subfigure}
\caption{Box plot of the simulated temperature deviation (with respect to the set point temperature) for all the houses, i.e., $\Delta_h$ values for all the houses in the set $\mathcal H$. The spread of the data  reveals two salient trends: (i) When using the ``min." objective function, simulation with a 1 hour look-ahead leads to temperatures with too big a spread in temperature deviation. However, by using a 3 hour look-ahead, one can achieve peak reduction of 5\%, 10\% or even 15\% (corresponding to curtailment factor values of $0.95, 0.90, 0.85$ respectively) with a deviation performance at least as good as the baseline model. (ii) When using the ``min-max" objective function, simulations with a 1 hour look-ahead lead to unacceptably large temperature deviations even though the spread is limited.  However, with a look-ahead of 3 hours, even 15\% peak reduction  is feasible with temperature deviation values and spread that are  much lesser than the baseline.}
    \label{fig:box}
\end{figure}

Consistent with results in Table \ref{tab:cen-times}, the box plots (Figure~\ref{fig:box}) also show that the ``min.'' model produces a deviation performance much worse than the ``min. max.'' model within the computational time limit imposed. The spread of the data in  Figure~\ref{fig:box} reveals two salient trends: (i) When using the ``min." objective function, a 1 hour look-ahead leads to too big a spread in temperature deviation data. However, by using a 3 hour look-ahead, one can achieve peak reduction of 5\%, 10\% or even 15\% (corresponding to curtailment factor values of $0.95, 0.90$, and $0.85$, respectively) $\gamma$ with a deviation performance at least as good as the baseline model. (ii) When using the ``min-max" objective function, a 1 hour look-ahead leads to unacceptably large temperature deviations even though the spread is limited.  However, with a look-ahead of 3 hours, even 15\% peak reduction  is feasible with   temperature deviation values and spread that are  much lesser than the baseline. This leads to a conclusion that if the technology and infrastructure is indeed set up for implementing the centralized models, it can provide the operator a more efficient way of providing equitable services to the customers without incurring a substantial increase in the natural gas consumption. Going further, the results from the centralized model allow one to conclude that equitable service can be provided while simultaneously reducing the peak natural gas consumption during  an extreme event scenario.

\section{Future Research Directions and Conclusion} \label{sec:future}
This paper develops an OCP for natural gas DR based on physics and thermodynamics, and demonstrates significant benefits of DR in the field of natural gas. In order to demonstrate the feasibility of a DR implementation for natural gas focused on residential loads (home heating), the residential loads are treated as the sole consumers of natural gas. The natural gas pipeline network is ignored, and two kinds of smart thermostats are modelled. 
Our formulation of the baseline OCP is reflective of the current state of practice, and our formulation has been compared  to the current  state for two different instances - a typical winter day where the maximum and minimum temperature attained  are  15 degrees  apart, and an extreme weather event represented by a polar vortex where the temperature range is almost three times as large. It has been shown that our decentralized formulation can lead to appreciable savings in gas consumption for a typical day, and do slightly better than current practice for days representative of winter storms. However, our centralized formulation can save up to 15\% gas consumption even on such days.

For both the centralized and decentralized look-ahead OCP models of DR, we have shown in the case studies that substantial reduction in gas savings for end customers can be achieved without compromising on the comfort level, while at the same time reducing peak consumption for the aggregator. The only cost to the customers is the thermostat installation in the case of the decentralized OCP and the communication network to communicate temperature inside the house to a central operator in the case of a centralized OCP. Hence, if this cost is offset using incentives, it will attract customers for this DR program. Thus, the aggregator could use the formulations to determine how much benefit they would get, and  then use this information to create an incentive that would get people to participate.  Currently, household customers have fixed rate plans or variable rate plans (whose prices change on a monthly basis)  so the only way to  lower the cost is to reduce consumption. Hence, our DR formulations are reflective of today's technologies and market structures, and more sophisticated price based DR will require the introduction of real-time gas rates.

To make the study more realistic in the  future, current research can be extended in the following  directions: (i) inclusion of network effects due to the pipeline network and the study of its impact on the proposed DR models, (ii) formulation and solution of  DR models that include uncertainty in the temperature forecast, (iii) data-driven estimation of the thermal coefficients of the house using measurements and the known uncertainty in the estimates of these coefficients, and finally (iv) quantify the link between tine-based rate/economic incentives and direct load control.

In conclusion, our article represents the start of an interesting and important problem of designing DR for a natural-gas only system. This work is timely given the realities of climate change; climate change is causing an increased number of  severe winter weather events and this in turn is causing  a stress on the natural gas pipeline networks, impeding the ability of these systems to provide equitable services to all its customers. The lack of resiliency of the natural-gas grid  results in inconvenience, disruption and even loss of lives. The proposed DR models are a simple technique to increase the resiliency of the natural gas systems, leading to more efficient operation in severe events as well as normal winter weather. Nevertheless, large scale implementation of these models requires major technological and infrastructure investment as well as policy changes that need to be initiated by the federal government.

\appendix
\section{Computation of thermal coefficients for a house} \label{sec:thermal-coeff}

We consider a house $h \in \mathcal H$ and let $\bm \theta_h^0$ denote the initial temperature inside house $h$. We let $\bm \Theta(t)$ denote the time-varying forecast of the ambient (outside) temperature. Each house has a gas furnace that burns fuel at the rate $\bm Q_h$  (\si{\kilogram\per\second}). Each house $h$ is also defined by $\bm V_h$ - the volume of the house  (\si{\cubic\metre}), $\bm \kappa_h$ - the heat transfer coefficient of the walls of the house based on insulation (\si{\watt\per\square\metre\per\kelvin}), $\bm A_h$ - the surface area of the walls of the house exposed to the outside temperature (\si{\square\metre}), and $\bar{\bm\theta}_h$ - the constant temperature set-point of the house (set by the occupants) (\si{\kelvin}). Finally, we let $\bm E^g$ (\si{\joule\per\kilogram}) denote the specific heating value of natural gas. When the furnace is on, the heat flux produced by the furnace is given by $\bm Q_h \cdot \bm E^g$ (\si{\joule\per\second}). Furthermore, the following statements are assumed about the physical system: the air inside the  house undergoes  an isochoric (constant volume) process, and hence the relevant specific heat of the gas $\bm C^a$  (\si{\joule\per\kilogram\per\kelvin}) is calculated at constant volume. We assume that $\bm C^a$ is constant over the temperature range of interest. We also neglect the thermal compressibility (density dependence on temperature) of air and assume a constant density $\bm \rho^a$ (\si{\kilogram\per\cubic\metre}) over the temperature range of interest. The only source of heat supply is the furnace, and heat is constantly lost through the walls of the house at a rate proportional to the surface area of the walls and the difference in  temperature between  the house and  the external environment.

Given these notations and assumptions, the dynamics of heat flow within the house $h$ is governed by:
\begin{flalign}
& \bm C^a \bm \rho^a \bm V_h \dot{\theta}_h(t) =  -\bm \kappa_h \bm A_h (\theta_h(t) - \bm \Theta(t)) + \bm Q_h \cdot \bm E^g \cdot q^g_h(t). \label{eq:ode} 
\end{flalign}
Comparing Eq. \eqref{eq:ode} and \eqref{eq:ode-h}, we can obtain expressions for the thermal coefficients $\bm \alpha_h$ and $\bm \beta_h$ as:
\begin{flalign}
\bm \alpha_h = \frac{\bm \kappa_h \bm A_h}{\bm C^a \bm \rho^a \bm V_h} \quad \text{ and }\quad \bm \beta_h = \frac{\bm E^g}{\bm C^a \bm \rho^a \bm V_h} \label{eq:thermal-coeff}
\end{flalign}
In all the computational experiments, the thermal coefficients for each house is computed using Eq. \eqref{eq:thermal-coeff}.

\section{List of abbreviations} \label{sec:abbreviations}
\vspace{2ex}
\begin{tabular}{ll}
    DR & Demand Response \\
    OCP & Optimal Control Problem \\
    FERC & Federal Energy Regulatory Committee \\
    ODE & Ordinary Differential Equation \\ 
    TCL & Thermostatically Controlled Load \\
    RH & Receding Horizon \\ 
    MILP & Mixed-Integer Linear Program \\ 
    IC & Initial Conditions \\
    TO & Timed Out
\end{tabular}

\section*{Acknowledgements}
This work was supported by the U.S. Department of Energy's Advanced Grid Modeling (AGM) project \emph{Joint Power System and Natural Gas Pipeline Optimal Expansion}. The research work conducted at Los Alamos National Laboratory is done under the auspices of the National Nuclear Security Administration of the U.S. Department of Energy under Contract No. 89233218CNA000001. We gratefully thank the AGM program manager Alireza Ghassemian for his support of this research.

\bibliographystyle{elsarticle-num-names}
\bibliography{references}

\end{document}